\documentclass[12pt]{amsart}
\usepackage{amsaddr}
\usepackage[margin=1in]{geometry}
\usepackage{wrapfig, caption}
\usepackage[utf8]{inputenc}
\usepackage{amsmath, graphicx}
\usepackage{color, bbold}
\usepackage{subcaption}
\usepackage{amsthm}
\usepackage{amssymb}
\usepackage{amscd}
\usepackage{mathtools}
\usepackage{hyperref}
\usepackage{tikz}
\usepackage{comment}
\usepackage{cleveref}
\usepackage{cancel}
\usepackage{graphicx}
\usepackage{microtype}
\usetikzlibrary{shapes.geometric}
\usepackage[shortlabels]{enumitem} % Sheila 2022 April 27 (warning: DO NOT USE with Beamer!)
\usepackage[numbers]{natbib}

\usetikzlibrary{patterns}
\tikzset{
  norm/.style     = {shape=circle, draw},
  blue/.style     = {shape=circle, draw, fill=blue!25},
  high/.style     = {shape=circle, draw, color=red},
  bluehigh/.style = {shape=circle, draw, color=red, fill=blue!25},
  red/.style      = {shape=circle, draw, fill=red!25},
  both/.style     = {shape=circle, draw, fill=violet!35},
  root/.style     = {node, bottom color=red!30},
  env/.style      = {treenode, font=\ttfamily\normalsize},
  dummy/.style    = {circle}
}
\tikzstyle{standard}=[circle, draw=black, fill=white, very thick, minimum size=7mm]
\tikzstyle{standard2}=[circle, draw=black, fill=white, very thick]
\tikzstyle{blue2}=[circle, draw=black, fill=blue!25, very thick]
\tikzstyle{small}=[circle, draw=black, fill=black, very thick, minimum size=4mm]
\tikzstyle{special}=[circle, draw=red!60, fill=red!5, very thick, minimum size=5mm]
\newtheorem{theorem}{Theorem}[section]
\newtheorem{lemma}[theorem]{Lemma}
\newtheorem{cor}[theorem]{Corollary}
\newtheorem{prop}[theorem]{Proposition}
\theoremstyle{definition}
\newtheorem{df}[theorem]{Definition}
\newtheorem{rem}[theorem]{Remark}

\title{Realizations and Uniqueness of Cut Complexes of Graphs}
\author[Y.F. Shen]{Yufeng Shen}
\address{Xi’an Jiaotong University, Xi’an, Shaanxi 710049, China, yufeng\_shen@stu.xjtu.edu.cn}
\author[Z.Y. Song]{Zhiyu Song}
\address{Nankai University, Tianjin 300071, China, 2210655@mail.nankai.edu.cn}
\author[F.L. Yu]{Fenglin Yu}
\address{Peking University, Beijing 100871, China, fenglin@stu.pku.edu.cn}
\author[W.H. Zhou]{Wuhan Zhou}
\address{Peking University, Beijing 100871, China, wuhanzhou@stu.pku.edu.cn}
\author[J.Q. Zhuang]{Jingqi Zhuang}
\address{Fudan University, Shanghai 200433, China, 22300680047@m.fudan.edu.cn}

\begin{document}
%\begin{abstract}
    %Cut complexes of graphs, defined as simplicial complexes whose facets are vertex subsets whose complements induce disconnected subgraphs, provide a bridge between graph theory and topological combinatorics. In this paper, we explore realizations and properties of these complexes.

    %First, we are inspired by a main theorem. For any $d$-dimensional simplicial complex on $n$ vertices, we construct a graph on $n + m(d,n)$ vertices such that the complex is the cut complex of the graph, and compute bounds on $m(d,n)$ for specific values of $d$ and $n$.  

    %The second part studies \emph{cut complexes} and how graph properties affect their structure. In particular, we establish a necessary and sufficient condition that characterizes the case when a $3$-cut complex uniquely determines the underlying graph.

    %Third, we present an algorithm to determine whether a given pure simplicial complex is the 3-cut complex of a special class of graphs.
    %\end{abstract}
\begin{abstract}
    In this paper, we investigate three fundamental problems regarding cut complexes of graphs: their realizability, the uniqueness of graph reconstruction from them, and their algorithmic recognition. We define the parameter \(m(d,n)\) as the minimum number of additional vertices needed to realize any \(d\)-dimensional simplicial complex on \(n\) vertices as a cut complex, and prove foundational bounds. Furthermore, we characterize precisely when a graph on \(n \ge 5\) vertices is uniquely reconstructible from its 3-cut complex. Based on this characterization, we develop an \(O(n^4)\) recognition algorithm. These results deepen the connection between graph structure and the topology of cut complexes.
\end{abstract}
\maketitle

\section{Introduction} \label{sec-1}
All graphs in this paper are finite and simple. We denote the vertex set and edge set of a graph \(G\) by \(V(G)\) and \(E(G)\), respectively.

A \emph{graph complex} is a simplicial complex built from a graph, encoding its combinatorial structure into a higher-dimensional topological object. Graph complexes appear in many branches of mathematics including commutative algebra (see \cite{commutative,Herzog2011,miller2004combinatorial}), representation theory (see \cite{bate2025}), and topological data analysis (see \cite{edelsbrunner2010computational,Holzinger2014}), reflecting their broad structural utility. In this paper, we focus on a recently introduced family of graph complexes called \emph{cut complexes}. In 2024, \cite{Bayer_2024} and \cite{Bayer_2024_02}, Bayer et al. introduced two new families of graph complexes called \emph{cut complexes} and \emph{total cut complexes}. Their work is motivated by a famous theorem of Ralf Fröberg \cite{Froberg_1990} connecting commutative algebra and graph theory through topology. %Our work further investigates several fundamental questions concerning their realizability, uniqueness, and algorithmic recognition.

\begin{figure}[htb]
\centering
\begin{subfigure}{0.4\textwidth}
\centering
\begin{tikzpicture}
%Nodes
\node[standard] (node1) at (18:1.3) {1};
\node[standard] (node2) at (90:1.3) {2};
\node[standard] (node3) at (162:1.3) {3};
\node[standard] (node4) at (234:1.3) {4};
\node[standard] (node5) at (306:1.3) {5};

%Lines
\draw (node1) -- (node4);
\draw (node2) -- (node4);
\draw (node2) -- (node5);
\draw (node5) -- (node3);
\draw (node3) -- (node1);
% \node at (0.2, -1) {The cycle $C_5$};
\end{tikzpicture}
\caption{$\Delta_{3}(C_5) = \langle 13,35,52,24,41 \rangle$}
\end{subfigure}
%$\Delta_2(C_5)$
%\begin{figure}[htb]
\begin{subfigure}{0.4\textwidth}
\centering
\begin{tikzpicture}
\draw[fill=gray!20] (0,0) -- (0.7,1.5) -- (1.4,0) -- (0,0);
%\draw[fill=gray!20] (2.1,1.5) -- (0.7,1.5) -- (1.4,0) -- (2.1,1.5);
%\draw[fill=gray!20] (2.1,1.5) -- (1.4,0) -- (2.8,0) -- (2.1,1.5);
\draw[fill=gray!20] (3.5,1.5) -- (2.8,0) -- (2.1,1.5) -- (3.5,1.5);
%\draw[fill=gray!20] (3.5,1.5) -- (4.2,0) -- (2.8,0) -- (3.5,1.5);
%\draw[very thick, color=red,->] (0,0) -- (.55*0.7,.55*1.5);
%\draw[very thick, color=red] (.55*0.7,.55*1.5) -- (0.7,1.5);
%\draw[very thick, color=red,->] (3.5,1.5) -- (3.5+.55*.7,1.5-.55*1.5);
%\draw[very thick, color=red] (3.5+.55*.7,1.5-.55*1.5) -- (4.2,0);

\node[standard] (n5) at (0,0) {5};
\node[standard] (n2) at (0.7,1.5) {1};
\node[standard] (n4) at (1.4,0) {3};
\node[standard] (n1) at (2.1,1.5) {6};
\node[standard] (n3) at (2.8,0) {4};
\node[standard] (n5') at (3.5,1.5) {2};
%
% \node at (2, -1) {$\Delta_2 (C_5)=\langle{245}{,124}{,134} {,135}{,235}\rangle$};
\end{tikzpicture}
\caption{$\Delta_3 ^t (C_6)=\langle{135}{,246}\rangle$}
\end{subfigure}
\caption{$\Delta_3(C_5)$ is $\mathbb{S}^1$ and $\Delta_{3} ^{t} (C_6)$ retracts to $\mathbb{S}^0$}
\label{fig:Delta2C5}
\end{figure}
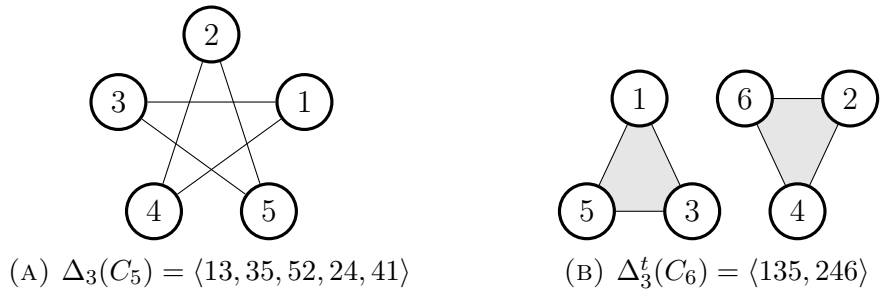

%Cut complexes of graphs serve as a fascinating intersection between graph theory and simplicial complexes in algebraic topology.
Given a graph $G$ on $n$ vertices with vertex set $V(G)$, the $k$-cut complex $\Delta_k(G)$ is defined as the simplicial complex on $V(G)$ whose facets are the $(n-k)$-subsets $S \subseteq V(G)$ such that the induced subgraph $G[V(G) \setminus S]$ is disconnected. Equivalently, the facets of $\Delta_k(G)$ are the complements of $k$-sets $U$ such that $G[V \setminus U]$ is disconnected, see Figure \ref{fig:Delta2C5} for an example.
These complexes capture topological properties of graphs, such as shellability and vertex decomposability, and have applications in understanding graph connectivity, reconstruction problems, and combinatorial topology.
A key result in this area establishes the universality of cut complexes for chordal graphs:
\begin{theorem} \label{con-thm} \cite{Bayer_2024}
    Let $\Delta$ be any pure simplicial complex. There exists some k and some chordal graph $G$ such that $\Delta = {\Delta}_{k}(G)$  
\end{theorem}
This theorem demonstrates that every pure simplicial complex can be realized as the cut complex of some chordal graph, highlighting the expressive power of chordal graphs in generating diverse topological structures. Building on this, we investigate three fundamental questions:
\begin{enumerate}
    \item \emph{Realization complexity:} How many extra vertices are needed to realize a given complex as a cut complex? We introduce the parameter \(m(d,n)\) to quantify this overhead.
    \item \emph{Uniqueness of reconstruction:} When does a cut complex determine the graph uniquely? We give a complete answer for \(3\)-cut complexes.
    \item \emph{Algorithmic recognition:} Can we efficiently decide whether a simplicial complex is a \(3\)-cut complex of some graph? We provide an \(O(n^4)\)-time algorithm.
\end{enumerate}

To formalize the first question, for a \(d\)-dimensional complex \(\Delta\) on \(n\) vertices, define \(m(d,n)\) as the minimal number such that there exists a graph \(G\) on at most \(n+m(d,n)\) vertices with \(\Delta = \Delta_{n+m(d,n)-(d+1)}(G)\). This parameter measures the vertex cost of embedding \(\Delta\) as a cut complex. We establish several bounds:
\begin{prop}
For any $n \geq 2$, $m(0,n) = 0$, $m(n-2,n) = 1$. Besides, for $n \geq 3$, $m(n-3,n) \leq \binom{n}{2} - 1$.
\end{prop}
\begin{prop}
For $n \geq 5$ and $1 \leq d \leq n-4$, we have $m(d,n) \geq 2$.
\end{prop}

We then study the uniqueness problem for \(3\)-cut complexes. Two vertices $u, v$ in a graph $G$ are \emph{twins} if they have identical open neighborhoods, i.e., $N(u) \setminus \{v\} = N(v) \setminus \{u\}$. A graph belongs to the family \(\mathcal{P}_4\) if it contains an induced path \(x-y-z-w\) with all vertices adjacent to a common set \(V'\). Our main uniqueness theorem states:

\begin{theorem}\label{szy}
Let $G$ be a graph on $n \geq 5$ vertices and $\Delta = \Delta _{3}(G)$, then for any graph on $n$ vertices, $G$ is the only graph that satisfies $\Delta_{3}(G) = \Delta$ if and only if there are no twins in $G$ and $G$ is not in $\mathcal{P}_4$.
\end{theorem}
This theorem implies that the presence of twins or membership in $\mathcal{P}_4$ allows for non-isomorphic graphs to share the same 3-cut complex, while their absence ensures uniqueness. The condition $n \geq 5$ avoids degenerate cases for smaller graphs where uniqueness may hold trivially or fail due to limited structure.

We also extend this to the total 3-cut complex $\Delta^t_3(G)$, defined similarly to $\Delta_3(G)$ but requiring that the induced subgraph on the complement of a face has an independent set of size $3$ (i.e., it is a total cut, where the induced subgraph of the complement has no edge), still see Figure \ref{fig:Delta2C5} for an example. A dominating pair in $G$ is a pair of vertices $\{u, v\}$ such that every other vertex in $G$ is adjacent to at least one of $u$ or $v$.
The analogous uniqueness result for total cut complexes is as follows:
\begin{theorem}\label{szy-total}
Let $G$ be a graph on $n \geq 3$ vertices and $\Delta = \Delta^t _{3}(G)$, then for any graph on $n$ vertices, $G$ is the only graph that satisfies $\Delta^t_{3}(G) = \Delta$ if and only if there is no dominating pair in $G$.
\end{theorem}
This characterization highlights the role of dominating pairs in obstructing uniqueness for total cut complexes, with a lower vertex threshold $n \geq 3$ reflecting the simpler conditions for smaller graphs in this variant. These theorems provide a complete classification for uniqueness in the $k=3$ case, building on the universality results for cut complexes of graphs. Proofs involve analyzing facet structures, neighborhood distinctions, and constructing counterexamples for non-unique cases.

Based on these characterizations, we develop an \(O(n^4)\) algorithm that decides whether a given pure simplicial complex is the \(3\)-cut complex of a graph without twins and not in \(\mathcal{P}_4\), and reconstructs the graph when possible.

The paper is organized as follows: Section~\ref{sec-2} studies \(m(d,n)\) in detail; Section~\ref{sec-3} proves the uniqueness theorems; Section~\ref{sec-4} presents the recognition algorithm; and the last remarks discuss open problems.

\section{definition and properties of \texorpdfstring{$m(d,n)$}{m(d,n)}} \label{sec-2}
General references for graph theory and simplicial complexes are \cite{TopMeth, bondy2011graph, Jonsson_2007}.

A natural question regarding cut complexes is how to determine whether a given pure simplicial complex can be realized as a (total) cut complex of a graph using a limited number of vertices. Based on Theorem~\ref{con-thm}, we start with the following basic definition.
\begin{df}
 For any $d$-dimensional complex $\Delta$ on $n$ vertices, define 
$\overline{\Delta} = \langle [n]\setminus F \mid F \text{ is a facet of } \Delta \rangle$, 
and $\overline{\Delta}^c = \langle F \in \tbinom{[n]}{n-d-1} \mid F \text{ is not a facet of } \overline{\Delta} \rangle$. Here, $\overline{\Delta}$ records the complements of facets of 
$\Delta$, and $\overline{\Delta}^c$ consists of all $(n-d-1)$-subsets that are not facets of $\overline{\Delta}$, intuitively, those subsets whose induced subgraphs would be connected if $\Delta$ is a cut complex.
\end{df}

The following two propositions help us to determine some necessary structural conditions that any cut complex must satisfy.
\begin{prop}\label{intersection1}
    Let $\Delta$ be a $d$-dimensional simplicial complex on $n$ vertices. If there is a graph $G$ with $n$ vertices such that $\Delta = {\Delta}_{n-(d+1)}(G)$, then for any two facets $A,B \in \overline{\Delta}^c$ with $|A \cap B| = 1$, there exists a sequence of facets $A=A_0,A_1,...,A_m=B \subset \overline{\Delta}^c$ satisfying $|A_i \cap A_{i+1}| = |A_i| - 1 = n-d-2$ for each $i$.
\end{prop}
\begin{proof}
    For a cut complex $\Delta = \Delta_{n-(d+1)}(G)$, the complex $\overline{\Delta}^c$ consists of all $n-(d+1)$-subsets $V \subset [n]$ that induce connected subgraph in $G$. Therefore, if we take $A, B \in \overline{\Delta}^c$, $|A \cap B| = 1$, then both $G[A]$ and $G[B]$ are connected and share exactly one vertex. Without loss of generality, suppose $G[A]$ and $G[B]$ are trees. Then we can choose a vertex $v \in G[A]$ and $u \in G[B]$ appropriately such that $G[(A\setminus \{v\}) \cup \{u\}]$ is connected. By repeating the operation above, we will get the required sequence.
\end{proof}
\begin{prop}\label{intersection0}
    Let $\Delta$ be a $d$-dimensional simplicial complex on $n$ vertices. If there is a graph $G$ with $n$ vertices such that $\Delta = {\Delta}_{n-(d+1)}(G)$, then for any two facets $A,B \in \overline{\Delta}^c$ with $|A \cap B| \neq \varnothing$, there exists a facet $C \subset \overline{\Delta}^c$ satisfying $|A \cap C| = |A| - 1 = n-d-2$.
\end{prop}
\begin{proof}
    As mentioned above, we can assume $G[A]$ and $G[B]$ are trees. Now let $X = A \setminus B$, $Y = B\setminus A$, we shall consider if $G[A]$ has leaves in $X$.
    
    First, if there is a vertex $v \in X$ such that $v$ is a leaf in $G[A]$, we can choose $u \in Y$ which is adjacent to $A \cap B$ and let $C = (A\setminus \{v\}) \cup \{u\}$. Then $C$ is what we want.
    
    Second, if $G[A]$'s leaves are all in $A \cap B$, then $|A \cap B| \geq 2$ and there is a leaf $v \in A \cap B$ such that there exists edges between $(A \cap B) \setminus \{v\}$ and $Y$. Now let $C = (A\setminus \{v\}) \cup \{u\}$, $u \in Y$ is adjacent to $(A \cap B) \setminus \{v\}$, $C$ is what we desired.
\end{proof}
\begin{prop}
Let $\Delta$ be a $d$-dimensional simplicial complex on $d+4$ vertices. Then $\Delta$ is the cut complex of a graph with $d+4$ vertices if and only if
\[
\left\langle F \in \binom{[d+4]}{d+1} \mid F \text{ is not a facet of } \Delta \right\rangle
\]
is the cut complex of a graph with $d+4$ vertices.
\end{prop}
\begin{proof}
    By noting that for a graph $G$ on $3$ vertices, $G$ is connected if and only if the complement of $G$, say $G^c$, is disconnected, we have $\Delta = \Delta_{3}(G)$ if and only if 
    \[
    \left\langle F \in \binom{[d+4]}{d+1} \mid F \text{ is not a facet of } \Delta \right\rangle = \Delta_{3}(G^{c})
    \]
\end{proof}
\begin{df}
    We write $m(d,n)$ to denote the minimal number such that for every $d$-dimensional simplicial complex $\Delta$ on $n$ vertices, there always exists a graph $G$ on no more than $n+m(d,n)$ vertices such that $\Delta = \Delta_{n+m(d,n) - (d+1)}(G)$. 
\end{df}
One can easily check that $m(d,n)$ exists for all $n \geq 2$ and $0 \leq d \leq n-2$. We have some basic propositions about $m(d,n)$.
\begin{prop}
    For any $n \geq 2$, $m(0,n) = 0$, $m(n-2,n) = 1$ and for $n \geq 3$, $m(n-3,n) \leq \binom{n}{2} - 1$.
\end{prop}
\begin{proof}
    For the first case, suppose $\Delta = \left\langle x_1,x_2,...,x_m\right\rangle$, $x_i \in [n]$, $1 \leq i \leq m$. If $m \leq n-2$, write $y_1,y_2,...,y_{n-m}$ be vertices that are not in $\Delta$, let graph $G$ has edges $E = \{ y_1x_1,x_1x_2,x_2x_3,...,x_{m-1}x_{m},x_{m}y_2,x_my_3,...,x_my_{n-m} \}$, then $\Delta = \Delta _{n-1}(G)$. (If $m = 1$ then we get a star, that is $E(G) = \{x_1y_1, x_1y_2,...,x_1y_{n-1}\}$) If $m = n-1$, let $y$ be the vertex that is not in $\Delta$, and $G$ with $G[\Delta] = K_{n-1}$ and no edge is adjacent to $y$, then $\Delta = \Delta _{n-1}(G)$. If $m = n$, then the empty graph on $n$ vertices is what we want.
    
    For the second case, note that for any graph on $n$ vertices, $\Delta _{1}(G) = \binom{[n]}{n-1}$, so $m(n-2,n) \geq 1$. On the other hand, let $G$ be a graph on $n+1$ vertices, then facets of $\Delta$ describe all the edges of $G$, hence $m(n-2,n) = 1$.
    The last case is a direct corollary of Theorem \ref{con-thm}. (In their proof, they construct a graph $G$ on $n + t$ vertices such that $\Delta _{n+t-(d+1)}(G) = \Delta$ for the $d$-dimensional $\Delta$ with $t$ facets.)
\end{proof}
A central lower-bound result is that for most dimensions, at least two extra vertices are necessary. The proof constructs explicit complexes that cannot be realized with only one extra vertex, leveraging the structural constraints given in Proposition~\ref{intersection1} and Proposition~\ref{intersection0}.
\begin{theorem}
\label{boundsfor_m(d,n)}
    For $n \geq 5$ and $1 \leq d \leq n-4$, we have $m(d,n) \geq 2$.
\end{theorem}
\begin{proof}
We construct a counterexample for each admissible $d$. First suppose $d\leq \frac{n-3}{2}$, then $n - (d+1) \geq d+2$.
Let $\Delta = \left\langle F \in \binom{[n]}{d+1} |F \neq \{1,2,...,d+1\}, \{n-d,n-d+1,...,n\} \right\rangle$. By using propositions above, we know that $\Delta$ cannot be a cut complex of a graph on $n$ vertices. We are going to prove $\Delta$ can't be a cut complex of a graph on $n+1$ vertices. Suppose $G$ is a graph on $n+1$ vertices, $V(G) = \{1,2,...,n,x\}$ and $\Delta = \Delta _{n-d}(G)$. Since $x$ is not the vertex in $\Delta$, we conclude that for any $A\in \binom{[n]}{n-d}$, $G[A]$ is connected. In other words $G[\{1,2,...,n\}]$ is $(d+1)$-connected. Hence, $d(v) \geq d+1$, $v \in [n]$. We use some claims to help us simplify the problem:
    
    \textbf{Claim}:
        For $G$ defined above, $x$ has at most one neighbor in $\{1,2,...,n-(d+1)\}$ or $\{d+2,d+3,...,n\}$
    \begin{proof}
        Because $G[\{1,2,...,n\}]$ is $d+1$-connected, $G[\{1,2,...,n-d\}]$ is connected. So there is a vertex $v \in \{1,2,...,n-(d+1)\}$ such that $G[\{1,2,...v-1,v+1,...,n-d\}]$ is still connected. Note that $G[\{1,2,...v-1,v+1,...,n-d,x\}]$ is not connected, so $x$ is not adjacent to $\{1,2,...v-1,v+1,...,n-d\}$. Hence we know that $x$ has at most one neighbor in $\{1,2,...,n-d\}$. The other situation is similar.
    \end{proof}
    \textbf{Claim}:
        For $G$ defined above, $x$ has exactly one neighbor in $\{d+2,d+3,...,n-(d+1)\}$
    \begin{proof}
        In fact, by using the connectivity of $G[\{1,2,...,n\}]$ we have $G[\{1,2,...,n-(d+1),u\}]$ is connected for $\forall u \in \{n-d,...,n\}$, hence $x$ cannot be adjacent to $\{n-d,...,n\}$. Similarly, $G[\{d+2,...,n,u\}]$ is connected for $\forall u \in \{1,...,d+1\}$, hence $x$ cannot be adjacent to $\{1,...,d+1\}$. To conclude, $x$ can only be adjacent to $\{d+2,d+3,...,n-(d+1)\}$ and have exactly one neighbor since $G[\{1,2,...,n-(d+1)\},x]$ is connected.
    \end{proof}
    \textbf{Claim}:
        For $n \geq 6$, $1 \leq d \leq \frac{n-3}{2}$, there is no graph on $n+1$ vertices such that $\Delta _{n-d}(G) = \left\langle F \in \binom{[n]}{d+1} |F \neq \{1,2,...,d+1\}, \{n-d,n-d+1,...,n\} \right\rangle$.
    \begin{proof}
        Otherwise, from the claim above, we may assume there is a graph $G$ on $\{1,2,...,n,x\}$ with $\Delta _{n-d}(G) = \left\langle F \in \binom{[n]}{d+1} |F \neq \{1,2,...,d+1\}, \{n-d,n-d+1,...,n\} \right\rangle$ and $x$ is adjacent to $d+2$. If $d \geq 2$, then $G[\{d,d+1,...,n-1,x\}]$ is connected since $G[\{1,2,...,n\}]$ is $d+1$-connected and $x$ is adjacent to $d+2$. However, there is a vertex $v \neq d+2$ in $\{d,d+1,...,n-1\}$ such that $G[\{d,d+1,...,v-1,v+1,...,n-1,x\}]$ is connected, and that is a contradiction since $\{d,d+1,...,v-1,v+1,...,n-1\} \neq \{1,2,...,n-(d+1)\},\{d+2,d+3,...,n\}$.
        
        If $d = 1$, then $G[\{1,2,3,5,6,...,n,x\}]$ is connected and there is a vertex $v \neq 3$ in $\{1,2,3,5,...,n\}$ such that $G[\{1,2,3,5,6,...,n,x\} \setminus \{v\}]$ is connected. But there is a contradiction since $\{1,2,3,5,6,...,n,\} \setminus \{v\} \neq \{1,2,...,n-2\},\{3,4,...,n\}$.
    \end{proof}
    \textbf{Claim}:
        For $n = 5$, $d = 1$, there is no graph on $6$ vertices such that 
        \[
        \Delta _{4}(G) = \left\langle\left. F \in \binom{[5]}{2} \right|F \neq \{1,2\}, \{4,5\} \right\rangle.
        \]
    \begin{proof}
        Otherwise, let $G$ on $\{1,2,3,4,5,x\}$ be the graph with 
        \[\Delta _{4}(G) = \left\langle\left. F \in \binom{[5]}{2} \right|F \neq \{1,2\}, \{4,5\} \right\rangle.
        \]
        From the claims above, we have $x$ is adjacent to $3$ and the degree of $3$ is at least 2 in $G[\{1,2,3,4,5\}]$. Because $G[\{1,2,3,x\}]$ and $G[\{3,4,5,x\}]$ are the only two connected induced subgraphs on $4$ vertices containing $x$, $3$ must have at least one neighbor in $\{1,2\}$ and $\{4,5\}$ respectively. Therefore, we have a connected induced subgraph $G[\{3,u,x,y\}]$, with $x\in \{1,2\}$, $y \in \{4,5\}$ and that is a contradiction.
    \end{proof}
    If $n \geq 6$, $\frac{n-2}{2} \leq d \leq n-4$, and let 
    $$\Delta = \left\langle F \in \binom{[n]}{d+1} | F \neq \{n-d,n-d+1,...,n\}, \{1,2,...,n-(d+2),2n-2d-1,...,n\}\right\rangle,$$
    we claim that there is no graph $G$ on $n+1$ vertices such that $\Delta = \Delta _{n-d}(G)$. First of all, the induced subgraph $G[\{1,2,...n\}]$ is $(d+1)$-connected and the Claims above still hold. Therefore, since $G[\{1,2,...,n-(d+1),x\}]$ and $G[\{n-(d+1),n-d,...,2n-2d-2,x\}]$ are connected, $x$ is adjacent to $\{1,2,...,n-(d+1)\} \cap \{n-(d+1),n-d,...,2n-2d-2\} = \{n-(d+1)\}$. However, $n-(d+1)$ has at least $d+1$ neighbors and $d+1 > n-(d+2)$, so $n-(d+1)$ is adjacent to $\{n-d,n-d+1,...,n\}$. Therefore, we have a connected induced subgraph $G[\{v_1,v_2,...,v_{n-d-3},n-(d+1),x,v\}]$, $\{v_1,v_2,...,v_{n-d-3},v\} \subset N(n-(d+1))$, $v \in \{n-d,n-d+1,...,n\}$, but that is a contradiction. ($N(v)$ means the neiborhood of $v$ in the given graph)
\end{proof}
%\begin{prop}
%    For any $n \geq 5$ and $1 \leq d \leq n-4$, we have $m(d,n) \leq m(d-1,n-1) + m(d,n-1)$.
%\end{prop}
\begin{prop}
    Let $\Delta$ be a $d$-dimensional simplicial complex and $V(\Delta) = [n]$, if $\Delta$ equals a total cut complex of a graph $G$ and $\bigcap \limits_{\text{F is a facet of }\Delta}F = \varnothing$, then there exists a graph $H$ on $n$ vertices such that $\Delta = \Delta_{n-(d+1)}^{t}(H)$.
\end{prop}
\begin{proof}
    If $G$ has $n$ vertices, then the statement is correct, so it suffices to prove the case where $G = (V, E)$ has more than $n$ vertices. We will prove that, for any vertex $v \notin [n]$, we have $v \nsim k$, $\forall k \in [n]$, and, hence, $v$ is an isolated vertex. Because $\bigcap \limits_{\text{F is a facet of }\Delta}F = \varnothing$, there exists a facet $F_{k}$ of $\Delta$ such that $k \notin F_{k}$. Therefore, $v,k \in V(\Delta) \setminus F_{k}$ and $v \nsim k$. Finally, we can check that $\Delta_{|V| - (d+1)}^{t}(G) = \Delta_{|V| - 1 - (d+1)}^{t}(G[V \setminus \{v\}])$, where $G[V \setminus \{v\}]$ is the induced graph of $G$ on $V \setminus \{v\}$.
\end{proof}
\section{uniqueness of 3 cut complex} \label{sec-3}
In this section, we study the uniqueness problem for \(3\)-cut complexes and then extend it to total $3$-cut complexes. We introduce the key concepts of \emph{twin vertices} and \emph{the family $\mathcal{P}_4$ of graphs}, proving several foundational propositions that illustrate how the presence of twins affects the uniqueness of the cut complex. 
\begin{df}
    For a graph $G$, we call two vertices $v$ and $u$ \emph{twins} if they have the same neighborhood or $N(u) \setminus \{v\} = N(v) \setminus \{u\}$
\end{df}
\begin{prop}\label{twin-prop-1}
    Let $G$ be a graph on $n \geq 4$ vertices, construct $G'$ by adding an edge $u \sim v$ to $G$, then $\Delta_{3}(G) = \Delta_{3}(G')$ \emph{if and only if} $u$ and $v$ are twins in $G$. 
\end{prop}
\begin{proof}
    Let $V$ denote the set of vertices of $G$ and $G'$.
    
    $\Rightarrow$: Otherwise, there is a vertex $w$ in $G$ such that $w$ is adjacent to one of $u$ and $v$. Without loss of generality, we assume that $w \sim u$, $w \nsim v$. Therefore, $V \setminus \{u,v,w\}$ is a facet of $\Delta_{3}(G)$ since $u \nsim v$ in $G$. However, $V \setminus \{u,v,w\}$ is not a facet of $\Delta_{3}(G')$ since $w \sim u$ and $u \sim v$ in $G'$, that is a contradiction.
    
    $\Leftarrow$: It suffices to prove $\Delta_{3}(G) \subset \Delta_{3}(G')$. Now suppose that $u$ and $v$ are twins in $G$ and $u \nsim v$, let $G'$ be the graph that adds $u \sim v$ on $G$. It is easy to see that we only need to consider facet $F \in \Delta _{3}(G)$ with $\{u,v\} \subset V \setminus F$. Define $S = \{F \text{ is facet of } \Delta_{3}(G) | \{u,v\} \subset V \setminus F\}$ and $A = \{x \in V |\exists F \in S,\text{ } V \setminus F =\{u,v,x\}\}$, because $u$ and $v$ are twins in $G$, $u$ and $v$ cannot be adjacent to $A$. Therefore, for $\forall x \in A$, $G'[\{u,v,x\}]$ is disconnected. To conclude, $V\setminus \{u,v,x\}$ is a facet of $\Delta_{3}(G')$ for $\forall x \in A$, so $S \subset \Delta_{3}(G')$ and $\Delta_{3}(G) \subset \Delta_{3}(G')$.
\end{proof}
Similarly, we have another proposition, the proof of which is similar to the above proposition.
\begin{prop}\label{twin-prop-2}
    Let $G$ be a graph on $n \geq 4$ vertices, construct $G'$ by removing an edge $u \sim v$ from $G$, then $\Delta_{3}(G) = \Delta_{3}(G')$ \textit{if and only if} $u$ and $v$ are twins in $G$. 
\end{prop}
\begin{cor}
    For graph $G$ and $H$ on $n \geq 4$ vertices, if $\Delta_{3}(G) = \Delta_{3}(H)$ and $E(G) \subsetneq E(H)$, then $G$ has twins.
\end{cor}

\begin{lemma} \label{key-lem}
    Let $G$ be a graph on $n \geq 4$ vertices and $V = V(G) = \{x,y,z,w,v_1,v_2,...,v_{n-4}\}$, we have following results:

    \begin{enumerate}
\item[(i)] If $V\setminus\{x,y,w\}, V\setminus\{x,y,z\}\in \Delta_{3}(G)$ and $V\setminus\{y,z,w\} \notin \Delta_{3}(G)$ or $V\setminus\{x,z,w\} \notin \Delta_{3}(G)$ then $x \nsim y$ in $G$.

\item[(ii)] If $V\setminus\{x,y,w\}, V\setminus\{x,y,z\} \notin \Delta_{3}(G)$ and $V\setminus\{y,z,w\} \in \Delta_{3}(G)$ or $V\setminus\{x,z,w\} \in \Delta_{3}(G)$ then $x \sim y$ in $G$.
\end{enumerate}
\end{lemma}
\begin{proof}
    We will only prove (i); the proof of (ii) is similar.

    Note that $\Delta_{3}(G) = \left\langle F \cap F' | F\text{, }F' \text{ are facets in } \Delta_{2}(G) \text{ and } |F \cap F'| = n-3\right\rangle$, if $x \sim y$ in $G$, then $V\setminus\{x,y\} \notin \Delta_{2}(G)$. 
    
    Therefore, $V\setminus\{y,w\}$, $V\setminus\{x,w\}$, $V\setminus\{y,z\}$, $V\setminus\{x,z\} \in \Delta_{2}(G)$ since $V\setminus\{x,y,w\}$, $V\setminus\{x,y,z\}$ $\in \Delta_{3}(G)$. 
    
    However, $V\setminus\{y,z,w\} = V\setminus\{y,w\} \cap V\setminus\{y,z\} \in \Delta_{3}(G)$ and $V\setminus\{x,z,w\} = V\setminus\{x,w\} \cap V\setminus\{x,z\} \in \Delta_{3}(G)$, which is a contradiction.
\end{proof}
\begin{df}
    Let $\mathcal{P}_4$ denote a set of graphs that have $4$ vertices $x,y,z,w$ and a subset $V'$ of $V(G)$ such that $N(x) = V' \cup \{y\}$, $N(y) = V' \cup \{x,z\}$, $N(z) = V' \cup \{y,w\}$, $N(w) = V' \cup \{z\}$.
\end{df}

\begin{figure}[ht]
    \centering
        \textbf{}
        \begin{tikzpicture}[scale=1, every node/.style={circle, draw, inner sep=2pt}]
            \node (A) at (180:1) {x};
            \node (B) at (120:1) {y};
            \node (C) at (60:1) {z};
            \node (D) at (1,0) {w};
            \node (E) at (0,-2) {V'};
            \foreach \i/\j in {A/B, B/C, C/D, A/E, B/E, C/E, D/E}
                \draw[thick] (\i) -- (\j);
        \end{tikzpicture}
        \caption{an example of a graph in $\mathcal{P}_{4}$}
\end{figure}
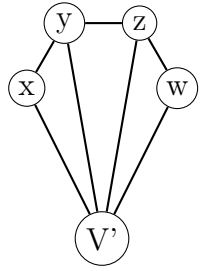

\begin{rem}
   When $V' = \varnothing$, the definition above means that $G$ has a connected component isomorphic to $P_4$, the path on the vertices of $4$. When $V \neq \varnothing$, we can see it as a generalization of fan graph, see \cite{2023arXiv230804512G}. 
\end{rem}

Now we state and prove the main result with respect to the uniqueness problem for $3$-cut complex.
\begin{theorem}\label{uniqueness}
    Let $G$ be a graph on $n \geq 5$ vertices and $\Delta = \Delta _{3}(G)$, then for any graph on $n$ vertices, $G$ is the only graph that satisfies $\Delta_{3}(G) = \Delta$ \emph{if and only if} there are no twins in $G$ and $G$ is not in $\mathcal{P}_4$.
\end{theorem}
\begin{proof}
    $\Rightarrow$: If there are twins $u$ and $v$ in $G$, then using Propositions \ref{twin-prop-1} and \ref{twin-prop-2} we can get another graph $G'$ such that $\Delta_{3}(G) = \Delta_{3}(G')$ by changing the relation between $u$ and $v$. Moreover, if $G$ is in $\mathcal{P}_4$, if $G$ has a connected component which is isomorphic to $P_{4}$, then from the figure $10$ we can construct another graph $G'$ such that $\Delta_{3}(G) = \Delta_{3}(G')$. The general case is similar to the above one.
    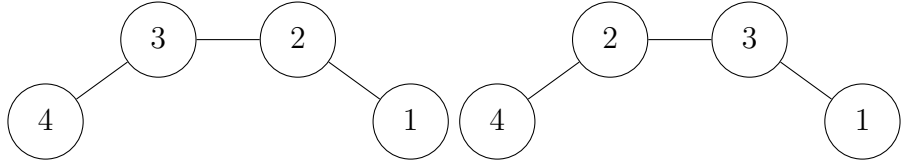
\begin{figure}[ht]
    \centering
    \begin{tikzpicture}[
    vertex/.style={draw, circle, minimum size=10mm, fill=white, 
                   text height=1.5ex, text depth=.25ex},
    every edge/.style={draw, thick}
    ]

    \foreach \i in {1,2,3,4} {
    \ifnum\i=1 \def\mylabel{$1$}\fi
    \ifnum\i=2 \def\mylabel{$2$}\fi
    \ifnum\i=3 \def\mylabel{$3$}\fi
    \ifnum\i=4 \def\mylabel{$4$}\fi
    
    \node[vertex] (v\i) at ({\i*36}:3cm) {\mylabel};
    }

    \foreach \i [evaluate={
    \next = int(mod(\i,4)+1);
    }] in {1,2,3} {
    \draw (v\i) -- (v\next);
    }

    \end{tikzpicture}
    \begin{tikzpicture}[
    vertex/.style={draw, circle, minimum size=10mm, fill=white, 
                   text height=1.5ex, text depth=.25ex},
    every edge/.style={draw, thick}
    ]

    \foreach \i in {1,2,3,4} {
        \ifnum\i=1 \def\mylabel{$1$}\fi
        \ifnum\i=2 \def\mylabel{$3$}\fi
        \ifnum\i=3 \def\mylabel{$2$}\fi
        \ifnum\i=4 \def\mylabel{$4$}\fi
    
        \node[vertex] (v\i) at ({\i*36}:3cm) {\mylabel};
    }

    \foreach \i [evaluate={
        \next = int(mod(\i,4)+1);
    }] in {1,2,3} {
        \draw (v\i) -- (v\next);
    }

    \end{tikzpicture}
    \caption{two paths on $4$ vertices that have the same $3$-cut complex}
    \end{figure}

    $\Leftarrow$: We will prove this  by verifying that every edge of $G$ is preserved when consider $\Delta_{3}(G)$. By noting that, and $G$ has no twins, so there is no other graph $H$ with $\Delta_{3}(G) = \Delta_{3}(H)$. Let $G$ be a graph that has no twins and none of the connected components of $G$ is isomorphic to $P_4$.

    \textbf{Claim}:
        For $G$ defined above and $u \sim v$, if $N(u) \cup N(v) = V(G)$, then for any $H$ with $\Delta _{3}(H) = \Delta_{3}(G)$ and $V(H) = V(G)$, $u \sim v$ in $H$.
    \begin{proof}
        Let $V = V(G)$ and $H$ is the graph mentioned above, if there are $x,y \in G$ such that $u \nsim x$ and $u \nsim y$, then we have $V \setminus \{u,v,x\}, V \setminus \{u,v,y\} \notin \Delta_{3}(H)$ and $V \setminus \{u,x,y\} \in \Delta_{3}(H)$, hence $u \sim v$ in $H$.

        If $u$ is adjacent to every vertex of $G$, then there is one vertex $x$ that is not adjacent to $v$ since $G$ has no twins. We may assume $x$ is the only vertex that is not adjacent to $v$, or it will be the case considered above. It is obvious that $x$ cannot be adjacent to every vertex of $N(v)$, so there is a vertex $y \neq u,v$ such that $x \nsim y$. Therefore $V \setminus \{u,v,x\}, V \setminus \{u,v,y\} \notin \Delta_{3}(H)$ and $V \setminus \{v,x,y\} \in \Delta_{3}(H)$, hence $u \sim v$ in $H$.

        If $|N(u)| = n-2$ and $x \nsim u$, then $x \sim v$ and there exists $y \in N(u)$ such that $x \nsim y$. Again we have $V \setminus \{u,v,x\}, V \setminus \{u,v,y\} \notin \Delta_{3}(H)$ and $V \setminus \{u,x,y\} \in \Delta_{3}(H)$, hence $u \sim v$ in $H$. 

        The above three cases are sufficient to prove our claim.
    \end{proof}
    Therefore, for $u \sim v$ in $G$, if $u$ or $v$ is not adjacent to at least two vertices in $G[N(u) \cup N(v)]$, then $u$ should be adjacent to $v$ in $H$. It suffices to prove $u \sim v$ will be preserved in $H$ if $u$ and $v$ are not adjacent to at most one vertex in $G[N(u) \cup N(v)]$.

    \textbf{Claim}:
        For $G$ defined above and $u \sim v$, if $N(u) \cup N(v) \neq V(G)$ and $u$ or $v$ is adjacent to every vertex in $N(u) \cup N(v)$, then for any $H$ with $\Delta _{3}(H) = \Delta_{3}(G)$ and $V(H) = V(G)$, $u \sim v$ in $H$.
    \begin{proof}
        Suppose $v$ is adjacent to every vertex in $N(u) \cup N(v)$. It suffices to consider the case when $u$ is not adjacent to exactly one vertex in $N(u) \cup N(v)$.

        If $N(u) \cup N(v) = \{u,v,x\}$, then $v \sim x$, $u \nsim x$. Because $G$ has no twins, there exists $y \in N(x)$ such that $y \neq v$. Moreover, since $G$ has no connected component isomorphic to $P_{4}$ and has $n \geq 5$ vertices, we can find a vertex $z$ which is adjacent to $x$ or $y$. No matter $z \sim x$ or $z \sim y$, $V \setminus \{x,y,z\} \notin \Delta_{3}(G)$. Observing that $V \setminus \{u,x,y\}, V \setminus \{u,x,z\} \in \Delta_{3}(H)$ and $V \setminus \{x,y,z\} \notin \Delta_{3}(H)$, $u \nsim x$ in $H$. In addition, $V \setminus \{u,v,x\} \notin \Delta_{3}(H)$, so $u \sim v$, $v \sim x$ in $H$.

        If $|N(u) \cup N(v)| \geq 4$. Let $x$ be the vertex that $u \nsim x$, $v \sim x$, if there exists vertex $y$ in $N(u)$ such that $x \nsim y$, then $u \sim v$ in $H$. If $x$ is adjacent to every vertex in $N(u)$, then because $G$ has no twins, $\exists y \in V \setminus (N(u) \cup N(v))$ such that $x \sim y$. Moreover, if x is adjacent to more than $1$ vertex in $V \setminus (N(u) \cup N(v))$, say $x \sim y \in V \setminus (N(u) \cup N(v))$, $x \sim y' \in V \setminus (N(u) \cup N(v))$, then $V \setminus \{v,x,y\}, V \setminus \{x,y,y'\} \notin \Delta_{3}(H)$, $V \setminus \{v,y,y'\} \in \Delta_{3}(H)$, hence $x \sim y$ in $H$. Note that $V \setminus \{u,x,y\} \in \Delta_{3}(H)$, so $u \nsim x$ in $H$ and $u \sim v$ in $H$.
        
        Now let's assume $x$ has only one neighbor in $V \setminus (N(u) \cup N(v))$, consider the neighborhood of $y$. 
        
        If $\exists z \in N(u) \cap N(v)$, $y \nsim z$, then $V \setminus \{v,x,y\}, V \setminus \{x,y,z\} \notin \Delta_{3}(H)$, $V \setminus \{v,y,z\} \in \Delta_{3}(H)$, hence $x \sim y$ in $H$. Moreover, $V \setminus \{u,x,y\} \in \Delta_{3}(H)$, so $u \nsim x$ in $H$ and $u \sim v$ in $H$.

        If $y$ is adjacent to every vertex in $N(u) \cap N(v)$, then we have $w \in V\setminus (N(u) \cup N(v) \cup \{y\})$ such that $y \sim w$. Otherwise, $\{u,v,x,y\}$ should be a good set to make $G$ an element of $\mathcal{P}_4$. Therefore, $V \setminus \{x,y,w\}, V \setminus \{x,y,v\} \notin \Delta_{3}(H)$ and $V \setminus \{v,y,w\} \in \Delta_{3}(H)$, so $x \sim y$ in $H$. Hence $u \nsim x$, $u \sim v$ in $H$.
    \end{proof}

    Now, let us organize the results we obtained. We have proved that for $u \sim v$ in $G$, if $u$ or $v$ is adjacent to all vertices in $(N(u) \cup N(v)) \setminus \{ u,v \}$, then $u \sim v$ should be preserved in $H$ when $\Delta _{3}(H) = \Delta_{3}(G)$. And if $u$ or $v$ is not adjacent to at least two vertices in $(N(u) \cup N(v)) \setminus \{ u,v \}$, then $u \sim v$ should also be preserved in $H$ when $\Delta _{3}(H) = \Delta_{3}(G)$. Therefore, there is only one last case left to consider. That is, when $u$ and $v$ are not adjacent to one vertex in $(N(u) \cup N(v)) \setminus \{ u,v \}$. We will now prove our last claim.

    \textbf{Claim}:
        For $G$ defined above and $u \sim v$, if $N(u) \cup N(v) \neq V(G)$ and $u$ and $v$ are not adjacent to one vertex in $N(u) \cup N(v)$, then for any $H$ with $\Delta _{3}(H) = \Delta_{3}(G)$ and $V(H) = V(G)$, $u \sim v$ in $H$.
    \begin{proof}
        Assume that $u \nsim y \in N(u) \cup N(v)$ and $v \nsim x \in N(u) \cup N(v)$, so we have $u \sim x$ and $v \sim y$. According to the previous proof, we only need to consider when $x$ and $y$ are adjacent to every vertex in $(N(u) \cup N(v)) \setminus \{ u,v \}$, hence $x \sim y$. Let $z$ be one vertex that is in $V \setminus (N(u) \cup N(v))$ such that $x \sim z$, $z$ must exist since $G$ has no twins. Similarly, there is a vertex $w \in V \setminus (N(u) \cup N(v))$ such that $w \sim y$.

        Actually, since $V \setminus \{u,x,z\}, V \setminus \{x,y,z\} \notin \Delta_{3}(H)$ and $V \setminus \{u,y,z\} \in \Delta_{3}(H)$, $x \sim z$ in $H$. Moreover, $V \setminus \{v,x,z\} \in \Delta_{3}(H)$ tells us $x \nsim v$ in $H$. Therefore, $u \sim v$ in $H$ since $V \setminus \{u,v,x\} \notin \Delta_{3}(H)$.
    \end{proof}
    
    To conclude, we have proved that for a graph $G$ without twins and none of its connected components is isomorphic to $P_4$, all its edges should be preserved in $H$ where $V(H) = V(G)$ and $\Delta_{3}(H) = \Delta_{3}(G)$. Finally, using the corollary $4.13$, we get the desired result that $E(H) = E(G)$.
\end{proof}

Moreover, we have a similar result on the $3$-total cut complex. Firstly, let us define dominating pairs of a graph.
\begin{df}
    For a graph $G$ and two vertices $u$ and $v$, we call them \emph{dominating pairs} if $N(u) \cup N(v) = V(G)$ or $N(u) \cup N(v) = V(G) \setminus \{u,v\}$.
\end{df}

Now we shall give the sufficient and necessary condition when the total $3$-cut complex is unique for graphs on $n$ vertices.

\begin{theorem}
    Let $G$ be a graph on $n \geq 3$ vertices and $\Delta = \Delta^t _{3}(G)$, then for any graph on $n$ vertices, $G$ is the only graph that satisfies $\Delta^t_{3}(G) = \Delta$ \emph{if and only if} there are no dominating pairs in $G$.
\end{theorem}
\begin{proof}
    $\Rightarrow$: Otherwise, if $G$ has dominating pair $u$ and $v$, then by changing the relation between $u$ and $v$ we get another graph $G'$. Without loss of generality, let $u \sim v$ in $G$ and $u \nsim v$ in $G'$. We only need to consider the $3$-set that contains $u$ and $v$. In fact, since $u$ and $v$ are dominating pair, for $\forall x \in V(G)$, $V(G) \setminus \{u,v,x\} \notin \Delta^t_{3}(G)$ and $V(G) \setminus \{u,v,x\} \notin \Delta^t_{3}(G')$. Therefore $\Delta^t_{3}(G) = \Delta^t_{3}(G')$ and that is a contradiction.
    
    $\Leftarrow$: Let $G$ be a graph  without dominating pairs, then for $\forall u, v \in G$, $\exists x \in G$ such that $u \nsim x$ and $v \nsim x$. If there is another graph $G'$ with $V(G) = V(G')$ and $\Delta^t_{3}(G) = \Delta^t_{3}(G')$, then for $u \nsim v$ in $G$ we claim that $u \nsim v$ in $G'$. Because $\exists x \in G$ such that $u \nsim v$, $v \nsim x$ and $x \nsim u$, $V(G) \setminus \{u,v,x\} \in \Delta^t_{3}(G)$. Hence $V(G) \setminus \{u,v,x\} \in \Delta^t_{3}(G')$ and $u \nsim v$. Moreover, for $u \sim v$ in $G$ we also have that $u \sim v$ in $G'$. Otherwise, if $u \nsim v$ in $G'$, then there must be a vertex $x \in V(G)$ such that $x \nsim u$ and $x \nsim v$ in $G$. Therefore $x \nsim u$ and $x \nsim v$ in $G'$, hence $V(G') \setminus \{u,v,x\} \in \Delta^t_{3}(G')$, which is contradict to $V(G) \setminus \{u,v,x\} \notin \Delta^t_{3}(G)$. Thus, we prove that $G = G'$.
\end{proof}

\section{algorithm} \label{sec-4}
Input a pure simplicial complex, and the following algorithm will determine whether it is a $3$-cut complex of a graph, suppose $\Delta$ has $n$ vertices and dimension $n-4$. We should input the connectivity of all induced graphs $G[V]$ where $V\in \binom{[n]}{3}$ and then the process has some steps.

Let us check if $u \sim v$ or $u \nsim v$ in $G$: 
\begin{itemize}
    \item [(i)] for $\forall x \in G$, check connectivity of $G[\{u,v,x\}]$ and let $A = \{x \in G|$ $G[\{u,v,x\}]$ $ \text{ is connected}\}$ $B = \{x \in G|G[\{u,v,x\}] \text{ is disconnected}\}$. If $A = \varnothing$, then $u \sim v$, if $B = \varnothing$ then $u \nsim v$.
    \item [(ii)] For $x, x' \in B$, if $G[\{u,x,x'\}]$ or $G[\{v,x,x'\}]$ is disconnected, then $u \nsim v$. For $y, y' \in B$, if $G[\{u,y,y'\}]$ or $G[\{v,y,y'\}]$ is connected, then $u \nsim v$.
    \item [(iii)] for $x\in B$, $y \in A$, we focus on the pairs $(x,y)$ such that $G[\{u,x,y\}]$ and $G[\{v,x,y\}]$ have different connectivity.
    
    When $G[\{u,x,y\}]$ disconnected and $G[\{v,x,y\}]$ connected: If there are at least two different $y$'s, then $u \sim v$. If there are at least two different $x$'s, then $u \nsim v$.

    When $G[\{u,x,y\}]$ connected and $G[\{v,x,y\}]$ disconnected: If there are at least two different $y$'s, then $u \sim v$. If there are at least two different $x$'s, then $u \nsim v$.
    
    If for $\forall (x,y)$, $G[\{u,x,y\}]$ and $G[\{v,x,y\}]$ are both connected or disconnected, return false.
    
    \item [(iv)] If there is only one pair such that $G[\{u,x,y\}]$ and $G[\{v,x,y\}]$ have different connectivity, then consider $(x,y) \times A \cup B$. $\{x,y,y'\}$ connected implies $u \sim v$ and $\{x,y,x'\}$ disconnected implies $u \nsim v$. Otherwise, return false.

    If there is $(x_1,y_1)$ such that $G[\{u,x_1,y_1\}]$ disconnected, $G[\{v,x_1,y_1\}]$ connected, and $(x_2,y_2)$ such that $G[\{u,x_2,y_2\}]$ connected, $G[\{v,x_2,y_2\}]$ connected: if $y_1 = y_2$, $u \sim v$, if $x_1 = x_2$, $u \nsim v$. If $y_1 \neq y_2$ and $x_1 \neq x_2$: $\{v,y_1,x_2\}$ disconnected implies $u \sim v$, connected implies $u \sim v$; $\{u, x_1, y_2\}$ disconnected implies $u \sim v$, connected implies $u \nsim v$.

    \item [(v)] Check the value of $(u,v)$.
\end{itemize}

%We will give the code at the end of our paper. There is only one question: whether there exists a simplicial complex that can give a graph but $\Delta \neq \Delta_{3}(G)$?

Actually, the algorithm only ensures that the result graph contains no twins and is not gem graph, so we should examine whether $\Delta_{3}(G) = \Delta$ after we get the result. As a result, we have an algorithm with complexity $O(n^4)$

\begin{theorem}\label{thm4.1}
    The algorithm above will tell if a simplicial complex is $3$-cut complex of a graph without twins and not in $\mathcal{P}_4$. 
    
    If the complex is $3$-cut complex of a graph without twins and not in $\mathcal{P}_4$, then the algorithm will give the graph with complexity $O(n^4)$.
\end{theorem}
\begin{proof}
    We shall assume the complex is a $3$-cut complex of a graph $G =(E, V)$ without twins and not in $\mathcal{P}_4$ to help us decide the adjacency of each pair of vertices of $\Delta$. For $u,v \in V$, if $A = \varnothing$ but $u \nsim v$ in $G$, then $u$ and $v$ should be adjacent to any other vertex of $V$ hence $u$ and $v$ are twins. Similarly, if $B = \varnothing$ and $u \sim v$, then $u$ and $v$ should be non-adjacent to any vertex of $V$, $u$ and $v$ are both isolated vertices.
    
    For $A, B \neq \varnothing$, we can use Lemma \ref{key-lem} to know: if there are $x,x' \in B$ such that $G[\{u,x,x'\}]$ or $G[\{v,x,x'\}]$ is disconnected then $u \sim v$, and if there are $y,y' \in A$ such that $G[\{u,y,y'\}]$ or $G[\{v,y,y'\}]$ is connected then $u \nsim v$.

    When $A, B \neq \varnothing$ and for $x,x' \in B$, $y,y' \in A$, if we always have $G[\{u,x,x'\}]$ and $G[\{v,x,x'\}]$ connected and $G[\{u, y, y'\}]$ and $G[\{v, y, y'\}]$ disconnected, we would like to hope there is at least one pair $(y,x) \in A \times B$ such that $G[\{u,x,y\}]$ and $G[\{v,x,y\}]$ have different connectivity. Otherwise, $u$ and $v$ will be twins in $G$. If $u \nsim v$, we have $u \sim x$, $v \sim x$ for $x \in B$, $u$ and $v$ are adjacent to at most one vertex in $A$. If $v \sim y$, $y \in A$, then $y \nsim y'$ $y \neq y' \in A$. Because $G[\{v,x,y\}]$ is connected for $\forall x \in B$, $G[\{u,x,y\}]$ is connected for $\forall x \in B$, hence $y \sim x$, $\forall x \in B$ and $v$ and $y$ are twins, that is a contradiction! The case $u \sim v$ is similar.

    Now we only need to consider the case $A, B \neq \varnothing$ and for $x,x' \in B$, $y,y' \in A$, we always have $G[\{u,x,x'\}]$ and $G[\{v,x,x'\}]$ connected and $G[\{u, y, y'\}]$ and $G[\{v, y, y'\}]$ disconnected. Moreover, there exists $(y, x) \in A \times B$ such that $G[\{u,x,y\}]$ and $G[\{v,x,y\}]$ have different connectivity. Let $C = \{(x,y)|G[\{u,x,y\}] \text{ is connected } \text{ and } G[\{v,x,y\}] \text{ is disconnected}\}$ $D = \{(x,y)|G[\{u,x,y\}] \text{ is disconnected } \text{ and } G[\{v,x,y\}] \text{ is connected}\}$. 
    
    \textbf{Claim:} If there are $(x_1,y_1),(x_2,y_2) \in C$ or $(x_1,y_1),(x_2,y_2) \in D$, and $y_1 \neq y_2$ then $u \sim v$.
    \begin{proof}
        If $u \nsim v$, $(x_1,y_1),(x_2,y_2) \in C$ and $y_1 \neq y_2$, the only possibility is $u \sim y_1$, $u \sim y_2$ (in this case, $u$ and $v$ are both adjacent to $x_1$ and $x_2$, so $x_i \nsim y_i$ since $G[u,x_i,y_i]$ and $G[v,x_i,y_i]$ have different connectivity) which is contradict to $G[\{u,y, y'\}]$ is disconnected for $y, y' \in A$. Hence $u \sim v$.
        Similarly, if we have $(x_1,y_1),(x_2,y_2) \in D$, and $y_1 \neq y_2$, we know that $u \sim v$ as well.
    \end{proof}
    
    Similarly, we have another conclusion.
    
    \textbf{Claim:} If there are $(x_1,y_1),(x_2,y_2) \in C$ or $(x_1,y_1),(x_2,y_2) \in D$, and $x_1 \neq x_2$ then $u \nsim v$.

    We now only need to consider the following case. $A, B \neq \varnothing$ and for $x,x' \in B$, $y,y' \in A$, we always have $G[\{u,x,x'\}]$ and $G[\{v,x,x'\}]$ connected and $G[\{u, y, y'\}]$ and $G[\{v, y, y'\}]$ disconnected and $|C|$, $|D| \leq 1$.

    If $|C \cup D| = 1$, say $C \cup D = \{(x,y)\}$. Without loss of generality, suppose $G[\{u,x,y\}]$ is disconnected, then let us consider the connectivity of $G[\{x,y,z\}]$, $z \in A \cup B$. 
    
    \textbf{Claim:} If $\exists z \in A$ such that $G[\{x,y,z\}]$ connected, then $u \sim v$. If $\exists z \in B$ such that $G[\{x,y,z\}]$ disconnected, then $u \nsim v$.
    \begin{proof}
        We will prove the case when $\exists z \in A$ such that $G[\{x,y,z\}]$ is connected. Otherwise, if $u \nsim v$, then we have $u \sim x$, $v \sim x$ and $u \nsim y$, $v \sim y$. Therefore $y \nsim x$ since $G[\{u,x,y\}]$ is disconnected, and $y \sim z$ and $x \sim z$ since $G[\{x,y,z\}]$ is connected. However, that will lead to a connected induced graph $G[\{v, y, z\}]$, $y, z \in A$, which contradicts our assumption.
    \end{proof}

    Still assume that $C = \varnothing$, $D = \{(x,y)\}$, if $G[\{x,y,z\}]$ is disconnected for $\forall z \in A$ and $G[\{x,y,z\}]$ is connected for $\forall z \in B$, then we have conclusion as follows.

    \textbf{Claim:} If there is $z \in A$ such that $G[\{u,x,z\}]$ is connected, then $u \nsim v$. Otherwise, if $G[\{u,x,z\}]$ is disconnected for $\forall z \in A$, then $u \sim v$.

    \begin{proof}
        If there is $z \in A$ such that $G[\{u,x,z\}]$ is connected. Noting that $G[\{u,x,y\}]$ and $G[\{x,y,z\}]$ are disconnected, we apply the lemma \ref{key-lem} to get $x \nsim y$. Therefore, $v \sim y$ and $v \sim x$, $u \nsim v$ since $G[\{u,v,y\}]$ is disconnected.

        If $G[\{u,x,z\}]$ is always disconnected when $z \in A$, we will show $u \sim v$ by contradiction. Assume $u \nsim v$, then $u \sim x$, $x \sim v$ and $v \sim y$ to ensure $D = \{(x,y)\}$. Moreover, one can check that $N(u) = B$, $N(x) = (B\backslash \{x\}) \cup \{u,v\}$, $N(v) = B \cup \{y\}$ and $N(y) = (B \backslash \{x\}) \cup \{v\}$ which means that $G$ is a graph in $\mathcal{P}_4$, that is a contradiction!
    \end{proof}

    When $A = \{y\}$, we will turn to the connectivity of $G[\{v, y, z\}]$, and the conclusion is as follows.

    \textbf{Claim:} If there is $z \in B$ such that $G[\{v,y,z\}]$ is disconnected, then $u \sim v$. Otherwise, if $G[\{v,y,z\}]$ is connected for $\forall z \in B$, then $u \nsim v$.

    The proof is similar to the previous one, so we will skip that.

    If $C = \{(x_1, y_1)\}$ and $D =\{(x_2,y_2)\}$. We will complete the whole proof after discussing the following cases. If $x_1 = x_2$, or $y_1 = y_2$, or $x_1 \neq x_2$, $y_1 \neq y_2$. The result is: if $x_1 = x_2$, then $u \nsim v$. If $y_1 = y_2$, then $u \sim v$. If $x_1 \neq x_2$, $y_1 \neq y_2$, by checking the connectivity of $G[\{u, x_1, y_2\}]$ we have: $u \nsim v$ if $G[\{u, x_1, y_2\}]$ is connected; $u \sim v$ if $G[\{u, x_1, y_2\}]$ is disconnected.

    \begin{itemize}
        \item [(1)] If $x_1 = x_2$, then $u \sim v$.
    \end{itemize}
    \begin{proof}
        Otherwise, if $u \sim v$, then $u \nsim y_1$, $v \nsim y_2$. Because $G[\{u,x_1,y_1\}]$ and $G[\{v,x_1,y_2\}]$ are connected, we have $u \sim x_1$, $v \sim x_1$, $y_1 \sim x_1$ and $y_2 \sim x_1$, which is contradict to $G[\{u,x_1, y_2\}]$ is disconnected.
    \end{proof}
    Similarly, we have the following conclusion.
    \begin{itemize}
        \item [(2)] If $y_1 = y_2$, then $u \sim v$.
    \end{itemize}

    In the end, if we can determine the relation between $(u,v)$ with $x_1 \neq x_2$, $y_1 \neq y_2$, then we have done so. Indeed, we have the result as below.
    \begin{itemize}
        \item [(3)] If $x_1 \neq x_2$, $y_1 \neq y_2$, then $u \nsim v$ if $G[\{u, x_2, y_1\}]$ is connected; $u \sim v$ if $G[\{u, x_2, y_1\}]$ is disconnected.
    \end{itemize}
    \begin{proof}
        We will prove this by contradiction. When $G[\{u, x_2, y_1\}]$ is connected and $u \sim v$. There must be $u \nsim x_2$ since $G[\{u, x_2, y_2\}]$ is disconnected and $G[\{v, x_2, y_2\}]$ is connected, which means $v \sim x_2$, $y_2 \sim x_2$, hence $u \nsim x_2$. Moreover, $u \nsim y_1$ since $u \sim v$ and $G[\{u,v,y_1\}]$ is disconnected. Therefore, $G[\{u, x_2, y_1\}]$ cannot be connected, that is a contradiction!

        Similarly, we have another conclusion that $u \sim v$ if $G[\{u, x_2, y_1\}]$ is disconnected.
    \end{proof}
    To conclude, we can judge the adjacency of $(u,v)$ after checking $A, B$, connectivity of $G[\{u,x,x'\}]$, $G[\{v,x,x'\}]$ where $x,x'$ both in $A$ or $B$, connectivity of $G[\{u,x,y\}]$, $G[\{v,x,y\}]$ where $x \in B$, $y \in A$ and pairs $(y,x) \in A \times B$ such that $G[\{u,x,y\}]$ and $G[\{v,x,y\}]$ have different connectivity. And then, to make sure the edge $u \sim v$ (or the non-edge $u \nsim v$) will not lead to contradiction, we should examine the connectivity of $G[\{u,v,x\}]$ with $x \neq u,v$ and $x \in V$. We can easily check the connectivity of a $3$-vertex induced graph by noting that it is connected if, and only if, it has at least $2$ edges.

    In the end, let us compute the complexity of the algorithm. Step $\text{(i)}$  has complexity $O(n)$, step $\text{(ii)}$ and step $\text{(iii)}$ have complexity $O(n^2)$ in total and step $\text{(iv)}$ has complexity of $O(n)$. The algorithm of checking the connectivity of $G[\{u,v,x\}]$ has complexity of $O(n)$. Therefore, combining a graph has $\binom{n}{2}$ pairs to be judged, the total complexity is $[O(n) + O(n^2) +O(n) + O(n)] \times O(n^2)= O(n^4)$.
\end{proof}

\section{remark} \label{sec-5}
In this paper, we define the parameter $m(d,n)$ to describe the minimum number of vertices of a graph, such that any $d$-dimensional simplicial complex $\Delta$ on $n$ vertices $[n]$ is the cut complex of a graph $G$ on $n+m(d,n)$ vertices, i.e. $\Delta = \Delta_{n+m(d,n)-(d+1)}(G)$. We give some basic properties of $m(d,n)$ by constructing specific $\Delta$ and give the exact value of $m(0,n)$ and $m(n-2,n)$.

To take a step forward, one may try to optimize the upper bound of $m(d,n)$ with fixed $d$; one upper bound of $m(d,n)$ is $\binom{n}{d+1} = O(n^{d+1})$ by using Theorem \ref{con-thm}.

As for $3$-cut complex of graphs, beyond our proof of the uniqueness of $3$-cut complex of a certain family of graphs, there are more questions about the construction of $3$-cut complex and $n$-cut complex with $n \geq 4$.

\textbf{Question:} Suppose we have two operations on a graph $G$:
\begin{itemize}
    \item (a) Change the relation between $u,v \in G$ when $u$ and $v$ are twins in $G$.
    \item (b) If there are four vertices $x,y,z,w \in G$ and $V' \subset V(G)$ such that $N(x) = V' \cup \{y\}$, $N(y) = V' \cup \{x,z\}$, $N(z) = V' \cup \{y,w\}$, $N(w) = V' \cup \{z\}$, then remove edges $xy$, $zw$ and add edges $xz$, $yw$.
\end{itemize}
For any two graphs $G$ and $H$ on $n$ vertices $[n]$, is it correct that if $\Delta_{3}(G) = \Delta_{3}(H)$ then $G$ can be made $H$ after a series of operations (a) and (b)?

\textbf{Question:} Does there exist any similar consequence of the uniqueness of the $n$-cut complex for $n \geq 4$?
\vspace{4em}
\renewcommand*{\bibfont}{\footnotesize}
    \bibliography{references}

@article{Bayer_2024,
  author = {Bayer, Margaret and Denker, Mark and Milutinovi\'{c}, Marija Jeli\'{c} and Rowlands, Rowan and Sundaram, Sheila and Xue, Lei},
  title = {Topology of Cut Complexes of Graphs},
  journal = {SIAM Journal on Discrete Mathematics},
  volume = {38},
  number = {2},
  pages = {1630-1675},
  year = {2024},
  doi = {10.1137/23M1569034},
  URL = {https://doi.org/10.1137/23M1569034},
  eprint = {https://doi.org/10.1137/23M1569034}
}

@article{Bayer_2024_02, 
    title={Total Cut Complexes of Graphs}, 
    volume={73}, 
    url={http://dx.doi.org/10.1007/s00454-024-00630-4}, 
    DOI={10.1007/s00454-024-00630-4}, 
    number={2}, 
    journal={Discrete \& Computational Geometry}, 
    publisher={Springer Science and Business Media LLC}, 
    author={Bayer, Margaret and Denker, Mark and Milutinović, Marija Jelić and Rowlands, Rowan and Sundaram, Sheila and Xue, Lei}, 
    year={2024}, 
    month={Feb}, 
    pages={500–527}, 
    language={en} 
}

@article{Froberg_1990,
  author = {Ralf Fröberg},
  journal = {Banach Center Publications},
  language = {eng},
  number = {2},
  pages = {57-70},
  title = {On Stanley-Reisner rings},
  url = {http://eudml.org/doc/268115},
  volume = {26},
  year = {1990},
}

@ARTICLE{2023arXiv230804512G,
       author = {{Grinberg}, Darij},
        title = "{An introduction to graph theory}",
      journal = {arXiv e-prints},
         year = 2023,
        month = aug,
          eid = {arXiv:2308.04512},
        pages = {arXiv:2308.04512},
          doi = {10.48550/arXiv.2308.04512},
archivePrefix = {arXiv},
       eprint = {2308.04512},
 primaryClass = {math.HO},
       adsurl = {https://ui.adsabs.harvard.edu/abs/2023arXiv230804512G}
}

@article{commutative,
    author = {Dochtermann, Anton and Engström, Alexander},
    title = {Cellular resolutions of cointerval ideals},
    journal = {Mathematische Zeitschrift},
    year = {2012}
}

@Inbook{Herzog2011,
author="Herzog, J{\"u}rgen
and Hibi, Takayuki",
title="Monomial Ideals",
bookTitle="Monomial Ideals",
year="2011",
publisher="Springer London",
address="London",
pages="3--22",
isbn="978-0-85729-106-6",
doi="10.1007/978-0-85729-106-6_1",
url="https://doi.org/10.1007/978-0-85729-106-6_1"
}

@book{miller2004combinatorial,
  title={Combinatorial Commutative Algebra},
  author={Miller, E. and Sturmfels, B.},
  isbn={9780387223568},
  lccn={04052495},
  series={Graduate Texts in Mathematics},
  url={https://books.google.co.jp/books?id=CqEHpxbKgv8C},
  year={2004},
  publisher={Springer New York}
}

@misc{bate2025,
      title={Homology of matching complexes and representations of symmetric groups}, 
      author={Michael Bate and Brent Everitt and Sam Ford and Eric Ramos},
      year={2025},
      eprint={2312.13750},
      archivePrefix={arXiv},
      primaryClass={math.GR},
      url={https://arxiv.org/abs/2312.13750}, 
}

@Inbook{Holzinger2014,
author="Holzinger, Andreas",
title="On Topological Data Mining",
bookTitle="Interactive Knowledge Discovery and Data Mining in Biomedical Informatics: State-of-the-Art and Future Challenges",
year="2014",
publisher="Springer Berlin Heidelberg",
address="Berlin, Heidelberg",
pages="331--356",
isbn="978-3-662-43968-5",
doi="10.1007/978-3-662-43968-5_19",
url="https://doi.org/10.1007/978-3-662-43968-5_19"
}

@book{edelsbrunner2010computational,
  title={Computational Topology: An Introduction},
  author={Edelsbrunner, H. and Harer, J.},
  isbn={9780821849255},
  lccn={2009028121},
  series={Applied Mathematics},
  url={https://books.google.co.jp/books?id=MDXa6gFRZuIC},
  year={2010},
  publisher={American Mathematical Society}
}

@book{Jonsson_2007,
  title={Simplicial Complexes of Graphs},
  author={Jonsson, J.},
  isbn={9783540758587},
  lccn={2007937408},
  series={Lecture Notes in Mathematics},
  url={https://books.google.com.hk/books?id=VCdHp6zlmKYC},
  year={2007},
  publisher={Springer Berlin, Heidelberg},
  volume = {1928}
}

@inbook{TopMeth,
  author = {Bj\"{o}rner, A.},
  title = {Topological methods},
  year = {1996},
  isbn = {0262071711},
  publisher = {MIT Press},
  address = {Cambridge, MA, USA},
  booktitle = {Handbook of Combinatorics (Vol. 2)},
  pages = {1819–1872},
  numpages = {54}
}

@book{bondy2011graph,
  title={Graph Theory},
  author={Bondy, A. and Murty, U.S.R.},
  isbn={9781846289699},
  lccn={2007923502},
  series={Graduate Texts in Mathematics},
  url={https://books.google.co.jp/books?id=HuDFMwZOwcsC},
  year={2011},
  publisher={Springer London}
}
    \bibliographystyle{apalike}
    % \tiny\bibliographystyle{alpha}
    
\end{document}